\documentclass[review]{elsarticle}
\usepackage{graphicx}
\usepackage{epstopdf}
\usepackage{lineno,hyperref}
 \usepackage{geometry}
\usepackage{epsfig}
\usepackage{lineno,hyperref}
\usepackage{bookmark}
\usepackage{multirow}
 \usepackage{amsfonts}
\usepackage{multirow}
 \usepackage{algorithm}
 \usepackage[figuresright]{rotating}
 \usepackage{amscd}
\usepackage{mathrsfs}
\usepackage{booktabs,longtable}
\usepackage{indentfirst}
\usepackage{subfigure}
\usepackage{amstext}
\usepackage{amssymb}
\usepackage{fancyhdr}
\usepackage{amsmath}
\modulolinenumbers[5]
\biboptions{sort&compress}
\newtheorem{thm}{Theorem}
\newtheorem{lem}[thm]{Lemma}
\newdefinition{rmk}{Remark}
\newproof{pf}{Proof}
\newproof{pot}{Proof of Theorem \ref{thm2}}
\newtheorem{alg}{Algorithm}
\newdefinition{exa}{Example}
\newtheorem{Def}{Definition}

\journal{Journal of \LaTeX\ Templates}

\begin{document}

\begin{frontmatter}

\title{ A Count Sketch  Kaczmarz Method For Solving Large Overdetermined Linear Systems\tnoteref{mytitlenote}}
\tnotetext[mytitlenote]{The work is supported by the National Natural Science Foundation of China (No. 11671060) and the Natural Science Foundation Project of CQ CSTC (No. cstc2019jcyj-msxmX0267)}

\author{Yanjun Zhang, \ Hanyu Li\corref{mycor}}
\cortext[mycor]{Corresponding author. E-mail addresses: yjzhang@cqu.edu.cn;\ lihy.hy@gmail.com or hyli@cqu.edu.cn.}

\address{College of Mathematics and Statistics, Chongqing University, Chongqing 401331, P.R. China}
\begin{abstract}
In this paper, combining  count sketch and maximal weighted residual Kaczmarz method, we propose a fast randomized algorithm for large overdetermined linear systems. Convergence analysis of the new algorithm is provided. Numerical experiments show that, for the same accuracy, our method behaves better in computing time compared with the state-of-the-art algorithm.
\end{abstract}


\end{frontmatter}

\linenumbers

\section{Introduction}
We consider the following consistent linear systems
\begin{equation}
\label{1}
Ax=b,
\end{equation}
where $A\in R^{m\times n}$ with $m\gg n $, $b\in R^{m}$, and $x$ is the $n$-dimensional unknown vector. As we know, the Kaczmarz method \cite{kaczmarz1} is a popular so-called row-action method for solving the systems (\ref{1}). In 2009, Strohmer and Vershynin \cite{Strohmer2009} proved the linear convergence of the randomized Kaczmarz (\texttt{RK}) method. 
Latter, many Kaczmarz type methods were proposed for different possible systems settings; see for example \cite{Needell2010, Eldar2011, Completion2013, Completion2015, Liu2019, Dukui2019, Wu2020, Chen2020} and references therein.

Recently, Bai and Wu \cite{Bai2018} constructed a greedy randomized Kaczmarz (\texttt{GRK}) method by introducing an efficient probability criterion for selecting the working rows from the coefficient matrix, which avoids a weakness of the one adopted in the \texttt{RK} method. 
Based on \texttt{GRK} method, a so-called relaxed greedy randomized Kaczmarz (\texttt{RGRK}) method was proposed in \cite{Bai2018r} by introducing a relaxation parameter, which makes the convergence factor of \texttt{RGRK} method be smaller than that of \texttt{GRK} method when the relaxation parameter $\theta\in[\frac{1}{2}, 1]$, and the convergence factor reaches the minimum  when $\theta=1$. For the latter case, i.e., $\theta=1$, Du and Gao \cite{Gao2019} called it the maximal weighted residual Kaczmarz (\texttt{MWRK}) method and carried out extensive experiments to test this method.

In this paper, inspired by dimensionality reduction techniques \cite{David}, we propose a count sketch Kaczmarz (\texttt{CSK}) method by combining count sketch \cite{Charikar2002, Thorup} and \texttt{MWRK} method. The convergence of \texttt{CSK} method is proved. Numerical experiments show that our method outperforms \texttt{MWRK} method in computing time.

The rest of this paper is organized as follows. In Section \ref{sec2}, some notations and the definition of count sketch are first given. Then, the \texttt{CSK} method is presented and its convergence is analyzed. Numerical experiments are given in Section \ref{sec3}.

\section{The \texttt{CSK} method }\label{sec2}
Throughout the paper, for a matrix $A$, $A^{(i)}$, $A_{(j)}$, $\sigma_i(A)$, $\sigma_r(A)$, $\|A\|_F$ and ${\rm R(A)}$ denote its $i$th row (or $i$th entry in the case of a vector), $j$th column, $i$th singular value, smallest nonzero singular value, Frobenius norm, and column space, respectively.

We now list the definition of count sketch which can be found in \cite{Charikar2002, Thorup}.

\begin{Def}{(Count Sketch transform).}
\label{countsketch}
 A count sketch transform is defined to be $\mathbf{S}=\Phi D\in{R^{d\times m}} $. Here, $D$ is an $m \times m$ random diagonal matrix with each diagonal entry independently chosen to be $+1$ or $-1$ with equal probability, and $\Phi\in \{0, 1\}^{d\times m} $ is a $ d\times m$ binary matrix with $\Phi_{h(i),i}=1$ and all remaining entries 0, where $h : [m] \rightarrow [d]$ is a random map such that for each $i\in [m]$, $h(i) = j$ with probability $1/d$ for each $j\in[d]$.
\end{Def}

Next, we give our new method.
\begin{alg}
\label{alg}
 The \texttt{CSK} method for the solution of the linear systems (\ref{1})
\begin{enumerate}[]
\item \mbox{INPUT:} ~ Matrix $A\in R^{m\times n}$, vector $b\in R^{m}$, parameter $d$, initial estimate $x_0$
\item \mbox{OUTPUT:} ~Approximate $x$ solving $Ax=b$
\item \mbox{Initialize:} ~Create a count sketch $\mathbf{S}\in{R^{d\times m}}$ with $d<m$, and compute $\widetilde{A}=\mathbf{S}A$ and $\widetilde{b}=\mathbf{S} b$.
\item For $k=0, 1, 2, \ldots,$ do until satisfy the stopping criteria
\item ~~~~Compute $i_{k}={\rm arg} \max \limits _{1 \leq i  \leq d}\{\frac{|\widetilde b^{(i)}-\widetilde A^{(i)} x_{k}| }{\| \widetilde A^{(i)}\|_{2}}\}$.
\item ~~~~Set $x_{k+1}=x_{k}+\frac{ \widetilde b^{ (i_{k}) }-\widetilde A^{ (i_{k} )} x_{k} }{\|\widetilde A^{(i_{k})}\|_{2}^{2}}(\widetilde A^{(i_{k})})^{T}$.
\item End for
\end{enumerate}
\end{alg}
\begin{rmk}
\label{rmkgl}
In \texttt{MWRK} method, the selection strategy for index is $i_{k}= \arg \max \limits _{1 \leq i  \leq m}\{\frac{|  b^{(i)}-  A^{(i)} x_{k}| }{\|   A^{(i)}\|_{2}}\}$. So, the difference between Algorithm \ref{alg} and \texttt{MWRK} method is that we introduce the count sketch transform $\mathbf{S}$. 
From \cite{David} and \cite{Clarkson2017}, 
we know that \textbf{S} can reduce the computation cost with keeping the most of the information of original problem. So, our method will behave better in runtime and a little worse in accuracy, which are conformed by numerical experiments given in Section \ref{sec3}. 
\end{rmk}

In the following, we provide theoretical guarantees for the convergence of the \texttt{CSK} method. A lemma is first given as follows, which plays a fundamental role in the convergence analysis.
\begin{lem}
\label{theorem1}(\cite{David})
If $\mathbf{S}$ is a count sketch transform with $O(n^2/(\delta\varepsilon^2))$ rows, where $0<\delta,\varepsilon<1$, then we have that
\begin{equation}
\label{lem11}
(1 -\varepsilon) \|A x\|^2_{2}\leq \|\mathbf{S} A x \|^2_{2}\leq(1 +
\varepsilon) \|A x\|^2_{2} ~for ~all~ x\in R^{n},
\end{equation}
 and
\begin{equation}
\label{lem12}
(1-\varepsilon) \sigma_{i}(\textbf{S}A)\leq \sigma_{i}(A) \leq(1 +\varepsilon)\sigma_{i}(\textbf{S}A) ~for ~all~ 1\leq i \leq n
\end{equation}
hold with probability $1-\delta$.
\end{lem}

\begin{thm}
\label{theorem2}
 Let $\mathbf{S}\in{R^{d\times m}}$ be a count sketch transform with $d=O(n^2/(\delta\varepsilon^2))$ and $x_\star=A^\dag b$ be the solution of the systems (\ref{1}). From an initial guess $x_0\in R^{ n}$ in the column space of $A^T$, for the sequence $\{x_k\}_{k=0}^\infty$ generated by the \texttt{CSK} method, we have that
\begin{eqnarray*}
\|x_{k+1}-x_\star\|^{2}_2\leq\left(1 -\frac{(1-\varepsilon)^3}{n}\cdot\frac{\sigma^2_{r}(A) }{\| A \|^2_2}\right)\| x_{k}-x_\star \|^2_2
\end{eqnarray*}
holds with probability at least $1-2\delta$.
\end{thm}

\begin{pf}
From Algorithm \ref{alg}, using the fact $Ax_\star=b$, we have
\begin{eqnarray*}
x_{k+1}-x_\star&=&x_{k}-x_\star+\frac{ \widetilde b^{(i_{k})}-\widetilde A^{(i_{k})} x_{k} }{\|\widetilde A^{(i_{k})}\|_{2}^{2}}(\widetilde A^{(i_{k})})^{T}
\\
&=& x_{k}-x_\star+\frac{ \mathbf{S}^{(i_{k})} b-\mathbf{S}^{(i_{k})}A  x_{k} }{\|\mathbf{S}^{(i_{k})}A\|_{2}^{2}}(\mathbf{S}^{(i_{k})}A)^{T}
\\
&=& \left(I- \frac{( \mathbf{S}^{(i_{k})}A)^T \mathbf{S}^{(i_{k})}A   }{\|\mathbf{S}^{(i_{k})}A\|_{2}^{2}}\right)( x_{k}-x_\star).
\end{eqnarray*}
Taking the square of the Euclidean norm on both sides and applying some algebra, we get
\begin{align}
\|x_{k+1}-x_\star\|^{2}_2
&=  \left\|\left(I- \frac{( \mathbf{S}^{(i_{k})}A)^T \mathbf{S}^{(i_{k})}A   }{\|\mathbf{S}^{(i_{k})}A\|_{2}^{2}}\right)( x_{k}-x_\star)\right \|^{2}_2 \notag
\\
&=( x_{k}-x_\star)^T \left (I- \frac{( \mathbf{S}^{(i_{k})}A)^T \mathbf{S}^{(i_{k})}A   }{\|\mathbf{S}^{(i_{k})}A\|_{2}^{2}}\right)( x_{k}-x_\star)\notag
\\
&=\|x_{k }-x_\star\|^{2}_2 -\frac{ | \mathbf{S}^{(i_{k})}A (x_{k }-x_\star)|^2  }{\|\mathbf{S}^{(i_{k})}A\|_{2}^{2}}. \label{eq:A}
\end{align}
Note that, from Algorithm \ref{alg},
\begin{eqnarray*}
i_{k}
&=&{\rm arg} \max \limits _{1 \leq i  \leq d}\{\frac{|\widetilde b^{(i)}-\widetilde A^{(i)} x_{k}| }{\| \widetilde A^{(i)}\|_{2}}\}
={\rm arg} \max \limits _{1 \leq i  \leq d}\frac{|  \textbf{S}^{(i)}b- \textbf{S}^{(i)}A x_{k}|^2 }{\|  \textbf{S}^{(i)}A\|^2_{2}}
\\
&=&{\rm arg} \max \limits _{1 \leq i  \leq d} \frac{ | \mathbf{S}^{(i)}A (x_{k }-x_\star)|^2  }{\|\mathbf{S}^{(i)}A\|_{2}^{2}}.
\end{eqnarray*}
Then
\begin{align}
\frac{ | \mathbf{S}^{(i_{k})}A (x_{k }-x_\star)|^2  }{\|\mathbf{S}^{(i_{k})}A\|_{2}^{2}}
&=\max \limits _{1 \leq i  \leq d} \frac{ | \mathbf{S}^{(i)}A (x_{k }-x_\star)|^2  }{\|\mathbf{S}^{(i)}A\|_{2}^{2}} \notag
 \geq \sum_{i=1}^{d}\frac{\|\mathbf{S}^{(i)}A\|^2_2}{\|\mathbf{S} A\|^2_F}\frac{ | \mathbf{S}^{(i)}A (x_{k }-x_\star)|^2  }{\|\mathbf{S}^{(i)}A\|_{2}^{2}} \notag
\\
&=\frac{ \|\mathbf{S} A(x_{k }-x_\star)\|^2_2}{\|\mathbf{S} A\|^2_F}.\label{eq:5}
\end{align}
Substituting  (\ref{eq:5}) into (\ref{eq:A}), we obtain
\begin{align}
\|x_{k+1}-x_\star\|^{2}_2
 &\leq \|x_{k }-x_\star\|^{2}_2 -\frac{ \|\mathbf{S} A(x_{k }-x_\star)\|^2_2}{\|\mathbf{S} A\|^2_F}
 \leq \|x_{k }-x_\star\|^{2}_2 -\frac{ \|\mathbf{S} A(x_{k }-x_\star)\|^2_2}{n\|\mathbf{S} A\|^2_2},\label{eq:6}
\end{align}
where the last inequality follows from the inequality $\|\mathbf{S} A\|^2_F\leq n\|\mathbf{S} A\|^2_2$.

As explained in  \cite{Bai2018}, since $x_\star=A^\dag b \in {\rm R(A^T)}$, by starting from an arbitrary initial guess $x_0$ in the column space of $A^T$, we have from the algorithm that $x_k$ also doe for each $k$ and hence, $x_k-x_\star$ is in the column space of $A^T$, which indicates that
$$\|  A(x_{k }-x_\star)\|^2_2\geq\sigma^2_{r}(A)\|x_{k }-x_\star\|^{2}_2.$$
Exploiting the above inequality and (\ref{lem11}), with probability $1-\delta$, we have
\begin{align}
\|\mathbf{S} A(x_{k }-x_\star)\|^2_2
 &\geq (1-\varepsilon)\|A(x_{k }-x_\star)\|^{2}_2\geq  (1-\varepsilon)\sigma^2_{r}(A)\|x_{k }-x_\star\|^{2}_2. \label{3}
\end{align}
Meanwhile, by (\ref{lem12}), with probability $1-\delta$, we have
\begin{eqnarray*}
 (1-\varepsilon) \sigma_{1}(\textbf{S}A) &\leq& \sigma_{1}(A).
\end{eqnarray*}
That is, with probability $1-\delta$, we have
\begin{align}
    \label{4}
 \|\textbf{S}A\|^2_2 &\leq \frac{1}{(1 -\varepsilon)^2} \| A \|^2_2.
\end{align}
Thus, combining (\ref{3}) and (\ref{4}), with probability at least $1-2\delta$, we get
\begin{align}
\frac{\|\mathbf{S} A(x_{k }-x_\star)\|^2_2}{\|\textbf{S}A\|^2_2 }
&\geq (1-\varepsilon)^3\cdot\frac{\sigma^2_{r}(A)\|x_{k }-x_\star\|^{2}_2  }{\| A \|^2_2}.\label{eq:7}
\end{align}
Substituting  (\ref{eq:7}) into (\ref{eq:6}), with probability at least $1-2\delta$, we have
\begin{eqnarray*}
\|x_{k+1}-x_\star\|^{2}_2
&\leq&\|x_{k }-x_\star\|^{2}_2 -\frac{(1-\varepsilon)^3}{n}\cdot\frac{\sigma^2_{r}(A)\|x_{k }-x_\star\|^{2}_2  }{\| A \|^2_2},
\end{eqnarray*}
which implies the desired result.
\end{pf}

\begin{rmk}
\label{rmk}
Note that $(1 -\frac{(1-\varepsilon)^3}{n}\cdot\frac{\sigma^2_{r}(A) }{\| A \|^2_2})<(1-\frac{\sigma^2_r(A) }{ \max \limits _{1\leq i \leq m} \sum\limits^{m}_{j=1,j\neq i}\|A^{(j)}\|^2_2 })$, where the latter is the convergence factor of \texttt{MWRK} method. So the convergence factor of \texttt{CSK} method is a little lager. This is because introducing count sketch transform $\mathbf{S}$ produces additional errors for algorithm. 
\end{rmk}

\section{Numerical experiments}\label{sec3}
In this section, we mainly compare the \texttt{CSK} method and the \texttt{MWRK} method in terms of the iteration numbers (denoted as ``IT'') and computing time in seconds (denoted as ``CPU''). We also report the iteration number speedup of \texttt{CSK} against \texttt{MWRK}, which is defined as $$\texttt{IT speedup}=\frac{\texttt{IT of }\texttt{MWRK}}{\texttt{IT of } \texttt{CSK}},$$ and the computing time speedup of \texttt{CSK} against \texttt{MWRK}, which is defined as $$\texttt{CPU speedup}=\frac{\texttt{CPU of }\texttt{MWRK}}{\texttt{CPU of } \texttt{CSK}}.$$

In all the following specific experiments, we generate the coefficient matrix $A\in R^{m\times n}$ and the solution vector $x_\star\in R^{n}$ using the MATLAB function \texttt{randn}, and the vector $b\in R^{m}$ by setting $b=Ax_\star$ and set $d=n^2$. We repeat 50 experiments and all the experiments start from an initial vector $x_0=0$, and terminate once the \emph{relative solution error }(RES), defined by $$\rm RES=\frac{\left\|x_{k}-x_\star\right\|^{2}_2}{\left\|x_\star\right\|^{2}_2},$$ satisfies $\rm RES\leq10^{-6}$, or the number of iteration steps exceeds 20000.


\begin{table}[tp]
  \centering
  \fontsize{6.5}{8}\selectfont
    \caption{ Numerical results for the \texttt{CSK} and \texttt{MWRK} methods.}
    \label{table1}
    \begin{tabular}{|c|c|c|c|c|c|c|}
    \hline
    \multirow{2}{*}{$m\times n$}&
    \multicolumn{3}{c|}{IT}&\multicolumn{3}{c|}{CPU}\cr\cline{2-7}
    &CSK&MWRK&IT speedup& CSK&MWRK&CPU speedup\cr
 \hline
  $300000\times 50$&     54.9000 &  31.0000 &   0.5647  &  0.2209 &   1.6878  &  7.6393 \cr\hline
 $300000\times 100$&     94.8600 &  63.0000  &  0.6641  &  0.4569 &   5.1097  & 11.1840\cr\hline
  $300000\times 150$&   132.7600 &  96.0000  &  0.7231  &  1.4894  & 10.3728 &   6.9645\cr\hline
  $400000\times 50$&     54.3600 &  29.0000  &  0.5335   & 0.2597  &  1.9737 &   7.6005 \cr\hline
 $400000\times 100$&        95.1400  & 60.0000  &  0.6306 &   0.5778   & 6.0903  & 10.5403\cr\hline
  $400000\times 150$&   132.5600 &  94.0000&    0.7091  &  1.8509  & 13.1016&    7.0783\cr\hline

 $500000\times 50$&     55.1000 &  29.0000  &  0.5263 &   0.3312   & 2.5878   & 7.8123 \cr\hline
 $500000\times 100$&      94.9600 &  60.0000   & 0.6318  &  0.7366  &  8.2903  & 11.2554\cr\hline
  $500000\times 150$&    132.5400  & 91.0000   & 0.6866  &  2.0513  & 17.3091  &  8.4383\cr\hline
      $600000\times 50$&    54.8800 &  28.0000 &   0.5102 &   0.4053   & 3.1672&    7.8142 \cr\hline
 $600000\times 100$&      95.1600 &  58.0000   & 0.6095  &  0.8566  & 10.2019 &  11.9103\cr\hline
  $600000\times 150$&   132.5800  &  92.0000  &  0.6939  &  2.3828 &  22.0388  &  9.2490\cr\hline
      $700000\times 50$&      54.7000 &  29.0000   & 0.5302  &  0.4550  &  3.9312  &  8.6401\cr\hline
 $700000\times 100$&     95.4200  & 58.0000  &  0.6078  &  1.0025   &12.2144  &\bf 12.1839\cr\hline
  $700000\times 150$&    132.3600 &  89.0000  &  0.6724  &  2.7194   &25.7228 &   9.4591\cr\hline

    \end{tabular}
\end{table}

The numerical results on IT and CPU are listed in Table \ref{table1}. Here, it should be pointed out that the IT and CPU in Table 1 denote the means of IT and CPU of 50 tests. From Table \ref{table1}, we see that the \texttt{CSK} method requires more iterations 
compared with the $\texttt{MWRK}$ method. This is because the \texttt{CSK} method has larger convergence factor and hence converges a little slower, which is consistent with the analysis of Remarks \ref{rmk} and \ref{rmkgl}. However, the runtime of the \texttt{CSK} method is less than that of the \texttt{MWRK} method, and the \texttt{CPU speedup} can be as large as 12.1839 in our experiments, 
which is consistent with the analysis of Remark \ref{rmkgl}.

\begin{figure}[ht]
 \begin{center}
\includegraphics [height=3.8cm,width=7.2cm ]{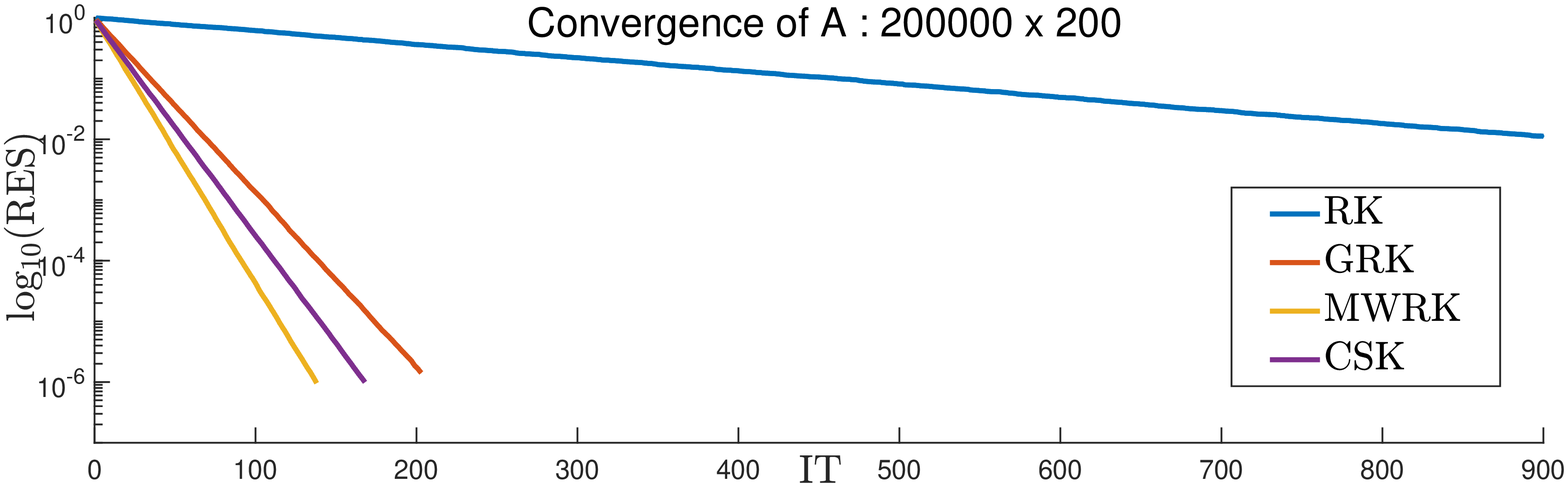}
 \includegraphics [height=3.8cm,width=7.2cm ]{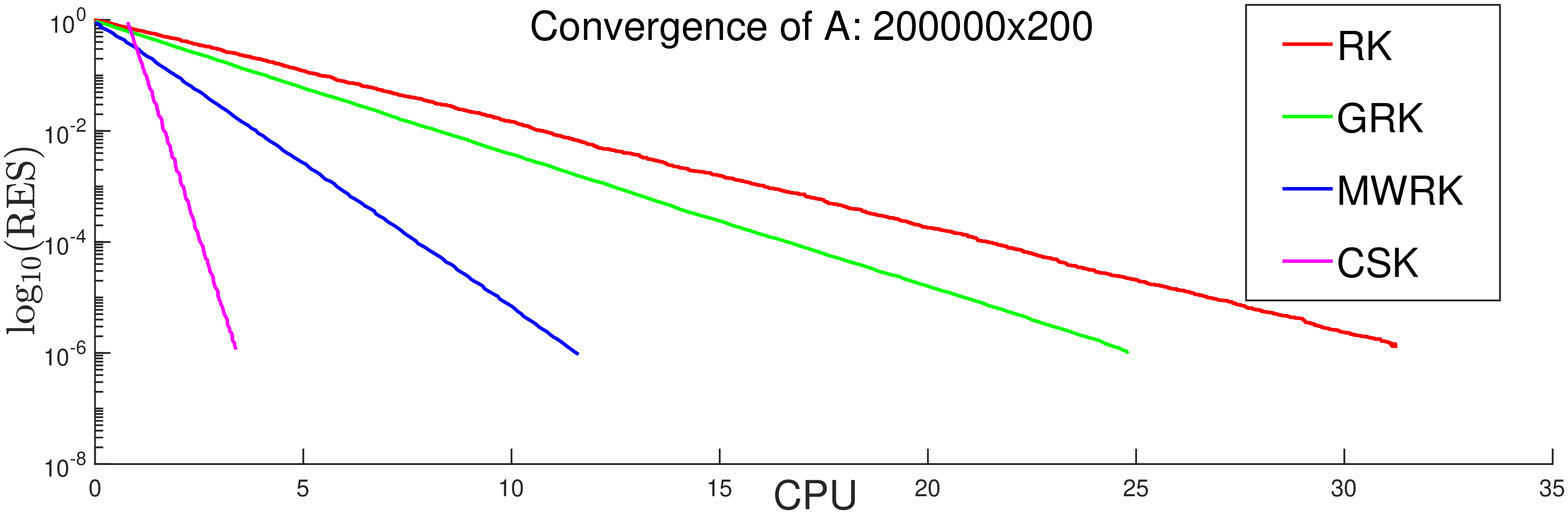}
 \end{center}
\caption{${\rm log_{10}(RES)}$ versus IT (left) and CPU (right) for \texttt{RK}, \texttt{GRK}, \texttt{MWRK}, and \texttt{CSK} when $A\in R^{200000\times200}$.}\label{fig}
\end{figure}

We also compare the performance of four algorithms (\texttt{RK}, \texttt{GRK}, $\texttt{MWRK}$, $\texttt{CSK}$). In Figure \ref{fig}, we plot the RES in base-10 logarithm versus the IT and CPU of four algorithms for $A\in R^{200000\times200}$. 
Each line represents the median RES at that iteration or CPU time over 50 trials. From the figure, we find that the $\texttt{CSK}$ and $\texttt{MWRK}$ methods outperform the \texttt{RK} and \texttt{GRK} methods in terms of IT and CPU, 
the $\texttt{MWRK}$ method converges fastest, 
and the $\texttt{CSK}$ method needs the least  runtime for the same accuracy. 



\bibliography{mybibfile}

\end{document}